\documentclass{article}

\usepackage{amssymb}
\usepackage{stmaryrd}
\usepackage[english,francais]{babel}

\newtheorem{theorem}{Theorem}[section]

\newtheorem{e-proposition}[theorem]{Proposition}

\newtheorem{e-definition}[theorem]{Definition\rm}

\newtheorem{theoreme}{Th\'eor\`eme}[section]

\newtheorem{proposition}[theoreme]{Proposition}

\def\Proof{\noindent \it Proof -- \rm}

\def\S{{\mathfrak  S}}

\def\<{\langle}
\def\>{\rangle}

\def\ashuff#1#2#3{
\kern 1pt \vrule height#1 \overline{\vrule height#3 width 0pt
\hskip#2} \rule{.3pt}{#1}\overline{\vrule height#3 width 0pt
\hskip#2} \rule{.3pt}{#1} \kern 1pt }

\def\Det{{\rm Det}}

\def\Pfa{{\rm Pf\,}}

\def\Det{{\rm Det}}

\def\binom#1#2{\left(#1\atop#2\right)}

\def\Carre3#1{\left[\begin{array}{ccc}#1\end{array}\right]}

\def\Pf#1{\mathtt{PF}^{(#1)}}

\def\og{\leavevmode\raise.3ex\hbox{$\scriptscriptstyle\langle\!\langle$~}}
\def\fg{\leavevmode\raise.3ex\hbox{~$\!\scriptscriptstyle\,\rangle\!\rangle$}}

\title{Hyperpfaffians}
\author{Ammar Aboud\footnote{USTHB, Faculty of Mathematics, Po. Box
32 El Alia 16111 Algiers, Algeria. aboudam@gmail.com}, Jean-Gabriel Luque\footnote{LITIS, Normandie Universit\'e,
Universit\'e de Rouen, Avenue de l'Universit\'e,
76801 Saint-\'Etienne du Rouvray Cedex,
France. jean-gabriel.luque@univ-rouen.fr}}
\begin{document}


\maketitle
\begin{abstract}
We define and inverstigate a generalization of the pfaffian for multiple array which interpolate between the hyperdeterminant and the hyperpfaffian.
\vskip 0.5\baselineskip
\end{abstract}


\selectlanguage{english}
 \section{Introduction}
 One of the simplest possible generalization of the determinant for higher-dimensional arrays is due to Cayley \cite{Cay1,Cay2} and 
 consists in considering a multiple alternating sum. The pfaffian of a skew-symmetric matrix is defined as the square root of the determinant. 
 In a more combinatorial way, it is also an signed sum but over perfect matchings instead of a signed sum over all the permutations. 
 A rather natural way  to define hyperpfaffian for $k$-tuple arrays
 consists in setting
 \begin{equation}
 \mathrm{HPf}(M):=\frac1{n!}\sum  \epsilon(\sigma)\prod_{i=1}^nM_{\sigma_1((i-1)k+1),\dots,\sigma_1(ik)},
 \end{equation}
 where the sum is over the permutations $\sigma\in\S_{nk}$ satisfying $\sigma(1)<\dots<\sigma(k)$, $\sigma(k+1)<\dots<\sigma(2k), \dots, \sigma((n-1)k+1)<\dots<\sigma(nk)$.
 This definition and a few variants are considered  in \cite{Abd,Barv,LT,Redelmeier}. In this paper we propose a more general definition for hyperpfaffian which interpolate
 between the hyperdeterminant and the hyperpfaffian. We prove several formulas (generalization Laplace expansion, hyperpfaffian of a sum, composition formula). Our main
 tool is the Grassmann-Berezin calculus. We also consider a generalization of a formula due to Gherardeli \cite{Ghe} relying hyperdeterminants and the Alon-Tarsi constant \cite{AT}.
 \section{Hyperpfaffians and Grassmann variables}
 \subsection{Combinatorial definition}
Let  $M=\left(M_{i_1,\dots,i_{mk}}\right)_{1\leq i_1,\dots,i_{mk}\leq mn}$ be a tensor. We define the following 
 polynomial which generalizes the combinatorial definition of the pfaffian of a matrix
 \begin{equation}\label{defpfaffian}
 \Pf m(M):=\frac1{n!}\sum  \epsilon(\sigma_1)\cdots\epsilon(\sigma_k)\prod_{i=1}^nM_{\sigma_1((i-1)m+1),\dots,\sigma_1(im),\dots,
 \sigma_k((i-1)m+1),\dots,\sigma_k(im)},
 \end{equation}
 where the sum runs over the $k$-tuples of permutations $(\sigma_1,\cdots,\sigma_k)\in \S_{mn}^k$ satisfying $\sigma_{j}((i-1)m+1)<\cdots<\sigma_j(im)$ for any $1\leq j\leq k$ and $1\leq i\leq n$.\\
 For any $M=\left(M_{i_1,\dots,i_{k}}\right)_{1\leq i_1,\dots,i_{k}\leq n}$, this definition allows to associate a polynomial $\Pf d(M)$ to $M$ for any $d$ which divides both $n$ and $k$.
 Notice that if $d=1$ we recover the  Cayley hyperdeterminant of $M$
 \begin{equation}
 \Pf1(M)=\Det(M)= \frac1{n!} \sum_{\sigma_1,\dots,\sigma_k\in\S_{n}}\epsilon(\sigma_1)\cdots\epsilon(\sigma_k)\prod_{i=1}^nM_{\sigma_1(i),\dots,\sigma_k(i)}.
 \end{equation}
 In the other end, if $k$ divides $n$ and $d=k$, we recover the notion of hyperpfaffian as defined in \cite{Barv,LT},
 \begin{equation}
 \Pf k(M)= \mathrm{HPf}(M)=\frac1{\left({n\over k}\right)!}\sum  \epsilon(\sigma)\prod_{i=1}^{n\over k}M_{\sigma((i-1)k+1),\dots,\sigma(ik)},
 \end{equation}
 where the sum runs over the permutations $\sigma\in\S_{n}$ satisfying $\sigma((i-1)k+1)<\dots<\sigma(ik)$ for any $1\leq i\leq {n\over k}$.
 Remark that if $d\neq n$ and $k$ is odd then $\Pf d(M)=0$.
  \subsection{Grassmann-Berezin calculus}
Let  $M=\left(M_{i_1,\dots,i_{mk}}\right)_{1\leq i_1,\dots,i_{mk}\leq mn}$ be a tensor where $mk$ is even and consider $k$ sets of formal
 variables $\eta^{(i)}=\{\eta^{(i)}_1,\dots,\eta^{(i)}_{mn}\}$ ($i=1\dots k$) satisfying the commutations $\eta^{(i)}_{j_1}\eta^{(i)}_{j_2}=-\eta^{(i)}_{j_2}\eta^{(i)}_{j_1}$ for any 
$1\leq i\leq k$, $1\leq j_1,j_2\leq mn$ and  $\eta^{(i_1)}_{j_1}\eta^{(i_2)}_{j_2}=\eta^{(i_2)}_{j_2}\eta^{(i_1)}_{j_1}$ for any 
$1\leq i_1\neq i_2\leq k$, $1\leq j_1,j_2\leq mn$. 
Let us introduce the notation known as Berezin integrals. The Berezin integral is a convenient tool for computing in Grasmann algebra (i.e., with anticommutative variables).
Let $f$ be a polynomials in the variables $\eta^{(1)},\dots,\eta^{(k)}$, we define
$
\int d\eta_{j_1}^{(i_1)}\cdots d\eta_{j_m}^{(i_m)}f:={\partial\over\partial \eta_{j_1}^{(i_1)}}\cdots
{\partial\over\partial \eta_{j_k}^{(i_k)}}f,
$
where each  $\partial\over\partial\eta_{j}^{(i)}$ acts on the Grasmann algebra as a left derivation ($\eta_{j}^{(i)}$ is pushed to the left, with a sign, and hence erased).
For simplicity we set also $\eta^{(i)}_J=\eta^{(i)}_{j_1}\cdots\eta^{(i)}_{j_m}$ for $J=\{{j_1}\leq\cdots\leq {j_m} \}$.\\
We define
$
\Omega_{m}(M):=\sum M_{i_1\dots i_{mk}} \eta^{(1)}_{\{i_1,\dots,i_m\}}\cdots \eta^{(k)}_{\{i_{(k-1)m+1},\dots,i_{km}\}},
$
where the sum is over the $km$-tuples $(i_1,\dots,i_{km})$ satisfying $i_1<\cdots<i_{m}, \dots,i_{(k-1)m+1}<\cdots<i_{km}$. By reorganizing 
the
monomials in the expansion of the polynomials one obtains the following result.
\begin{proposition}\label{Minors}
Let $\ell$ a divisor of $n$ and $I^{(1)},\dots, I^{(k)}$ be $k$ subsets of $\{1,\dots,mn\}$ of cardinality $\ell m$. One has
\begin{equation}
\frac1{\ell!}\int d\eta^{(1)}_{I^{(1)}}\cdots d\eta^{(k)}_{I^{(k)}}\Omega_m(M)^\ell=
\Pf m\left(M\left[\overbrace{I^{(1)}|\dots|I^{(1)}}^{\times m}|\cdots|\overbrace{I^{(k)}|\dots|I^{(k)}}^{\times m}\right]\right)
\end{equation}
where $M\left[I_1|\cdots|I_{mk}\right]$ is a hyperminor of $M$ that is the tensor obtained by selecting the entries
whose first index belongs in $I_1$, the second index belongs in $I_2$ etc. \\
As a special case, one has
\begin{equation}
\frac1{n!}\int d\eta^{(1)}_{\{1,\dots, mn\}}\cdots d\eta^{(k)}_{\{1,\dots,mn\}}\Omega_m(M)^n=
\int d\eta^{(1)}_{\{1,\dots, mn\}}\cdots d\eta^{(k)}_{\{1,\dots,mn\}}e^{\Omega_m(M)}
=\Pf m(M).
\end{equation}
\end{proposition}

 \section{Some formulas}
 \subsection{Generalization of the Lapace formula}
 Let $0<n'<n$ be an integer. We split $\Omega$ into two disjoint sums $\Omega_m(M)=\Omega'_m(M)+\Omega''_m(M)$
 where 
 \begin{equation}
\Omega'_{m}(M):=\sum M_{i_1\dots i_{mk}} \eta^{(1)}_{\{i_1,\dots,i_m\}}\cdots \eta^{(k)}_{\{i_{(k-1)m+1},\dots,i_{km}\}}
\end{equation}
where  the sum runs over the $km$-tuples $(i_1,\dots,i_{km})$ satisfying $i_1<\cdots<i_{m}, 
\dots,i_{(k-1)m+1}<\cdots<i_{km}$ and $i_1\in\{1,\dots,n'\}$.
The commutativity rules give
\begin{equation}
\Omega_m(M)^n=\binom n{n'}\Omega'_m(M)^{n'}\Omega''_m(M)^{n-n'}=\binom n{n'}\Omega'_m(M)^{n'}\hat\Omega''_m(M)^{n-n'},
\end{equation}
with
$
\hat\Omega''_{m}(M):=\sum M_{i_1\dots i_{mk}} \eta^{(1)}_{\{i_1,\dots,i_m\}}\cdots \eta^{(k)}_{\{i_{(k-1)m+1},\dots,i_{km}\}}
$
where the sum runs over the $km$-tuples $(i_1,\dots,i_{km})$ satisfying $i_1<\cdots<i_{m}, \dots,i_{(k-1)m+1}<\cdots<i_{km}$ and
 $i_1,\dots,i_m\not\in\{1,\dots,n'\}$.
But from Proposition \ref{Minors}, one has
\begin{equation}
\Omega'_m(M)^{n'}\hat\Omega''_m(M)^{n-n'}=n'!(n-n')!\sum^{\wedge}\Pf m\left(M[I]\right)\Pf m\left(M[J]\right)
\eta^{(1)}_{I^{(1)}}\cdots\eta^{(k)}_{I^{(k)}}\eta^{(1)}_{J^{(1)}}\cdots\eta^{(k)}_{J^{(k)}},
\end{equation}
where $\displaystyle\sum^{\wedge}$ means that the sum runs over the pairs of sequences 
 $$I=\left[\overbrace {I^{(1)}|\cdots|I^{(1)}}^{\times m}|\cdots|\overbrace {I^{(k)}|\cdots|I^{(k)}}^{\times m}\right]
\mbox{ and } J=\left[\overbrace {J^{(1)}|\cdots|J^{(1)}}^{\times m}|\cdots|\overbrace {J^{(k)}|\cdots|J^{(k)}}^{\times m}\right]$$
 satisfying $\{1,\dots,n'\}\subset I^{(1)}$, $\mathrm{card} (I^{(s)})=n'm$ and $J^{(s)}=\{1,\dots,nm\}\setminus I^{(s)}$ for any $1\leq s\leq k$.
 We deduce the following result which generalizes the Laplace expansion rule.
 \begin{theorem}\label{Laplace} One has
 \begin{equation}
 \Pf m(M)=\sum^{\wedge} (-1)^{\Lbag I\Rbag}\Pf m\left(M[I]\right)\Pf m\left(M[J]\right)
\end{equation}
where $\Lbag I\Rbag=k\binom{n'm+1}2+\sum_{i=1}^k\sum_{e\in I^{(i)}}e$.
 \end{theorem}
 \subsection{Hyperfpaffian of a sum}
 Let  $N=\left(N_{i_1,\dots,i_{mk}}\right)_{1\leq i_1,\dots,i_{mk}\leq mn}$ be another tensor. Since, $\Omega_m(M+N)=\Omega_m(M)+\Omega_n(M)$ one obtains
$
 \Omega_m(M+N)^n=\sum_{\ell=0}^n\binom n\ell\Omega_m(M)^\ell\Omega_m(N)^{n-\ell}$. Hence, proposition \ref{Minors} implies the following result.
 \begin{proposition} One has
 \begin{equation}
 \Pf m(M+N)=\sum_{\ell=0}^n\sum (-1)^{\Lbag I\Rbag}\Pf m(M[I])\Pf m(N[J]).
 \end{equation}
 where the second sum runs over the pairs of sequences
 $$I=\left[\overbrace {I^{(1)}|\cdots|I^{(1)}}^{\times m}|\cdots|\overbrace {I^{(k)}|\cdots|I^{(k)}}^{\times m}\right]
\mbox{ and } J=\left[\overbrace {J^{(1)}|\cdots|J^{(1)}}^{\times m}|\cdots|\overbrace {J^{(k)}|\cdots|J^{(k)}}^{\times m}\right]$$
 satisfying  $\mathrm{card} (I^{(s)})=\ell m$ and $J^{(s)}=\{1,\dots,nm\}\setminus I^{(s)}$ for any $1\leq s\leq k$. 
 \end{proposition}
 \subsection{Composition of Hyperpfaffians}
 Suppose now $m=pm'$, in this case $\Omega_{m'}(M)^{pn}=\left(\Omega_{m'}(M)^p\right)^n$. But from Proposition \ref{Minors}, one has
 \begin{equation}
 \Omega_{m'}(M)^p=p!\sum\Pf mM[\overbrace{I^{(1)}|\cdots|I^{(1)}}^{\times m'}|\cdots|\overbrace{I^{(k)}|\cdots|I^{(k)}}^{\times m'}]\eta_{I^{(1)}}^{(1)}\cdots
 \eta_{I^{(k)}}^{(k)},
 \end{equation}
 where the sum is over the sets $I^{(1)},\dots,I^{(k)}\subset\{1,\dots,nm\}$ of cardinality $m$. So $\Omega_{m'}(M)^p$ is written as $\Omega_{m}(M')$ where $M'$ 
 is a $mn^{\otimes m'k}$ tensor.
 More explicitely, applying again Proposition \ref{Minors}, we obtain 
 \begin{proposition} One has
 \begin{equation}
 \Pf {m}\left(\Pf {m'}(M\langle i_1,\dots,i_{mk}\rangle)\right)_{1\leq i_1,\dots,i_{mk}\leq mn}={\binom{np}{p,\cdots,p}\over n!}\Pf {m'}(M)
 \end{equation}
 with
 $$M\langle i_1,\dots,i_{mk}\rangle=M\left[\overbrace{\{i_1,\dots,i_m\}|\cdots|\{i_1,\dots,i_m\}}^{\times m'}|\cdots|
\overbrace{\{i_{m(k-1)+1},\dots,i_{mk}\}|\cdots|\{i_{m(k-1)+1},\dots,i_{mk}\}}^{\times m'} \right].$$
 \end{proposition}
 \section{Hyperpfaffians and generalized latin squares}
 A $(m,k)$-latin quasisquare is a $m\times mk$ matrix
 $$\left[
 \begin{array}{ccccccc}
 \sigma_1(1)&\cdots&\sigma_1(m)&\cdots& \sigma_k(1)&\cdots&\sigma_k(m)  \\
 \vdots      &    &\vdots            &\cdots&  \vdots&&\vdots\\
 \sigma_1((k-1)m+1)& \cdots &\sigma_1(km) &\cdots& \sigma_k((k-1)m+1)&  \cdots  &\sigma_m(km)        \\
 \end{array}
 \right] $$
 where each $\sigma_i$  is a permutation and each line is a permutation $\tau\in\S_{km}$ satisfying $\tau((\ell-1)m+1<\cdots<\tau(\ell m)$ for any $1\leq\ell\leq k$. We denote by $LQ(m,k)$ the set of the $(m,k)$-latin quasisquares.
 To each $c\in LQ(m,k)$ we associate a sign $\varepsilon(c)$ which is the product of the signs of the permutations $\sigma_i$'s and the signs of the lines.\\
 Let $\mathcal A^{n}$ be the unique $n^{\otimes n}$ antisymmetric tensor such that $\mathcal A^{n}_{1,\dots,n}=1$.  
 \begin{proposition}
 One has \begin{equation}
 \Pf m\left(\mathcal A^{mk}\right)=\frac1{k!}\sum_{c\in LQ(m,k)}\varepsilon(c)\end{equation}
 \end{proposition}
 \Proof
 Observe that
 \[
 \Omega_m\left(\mathcal A^{mk}\right)=\sum_{\tau\in\S_{mk}\atop
 \forall \ell\in\{1,\dots,k\}, \tau((\ell-1)m+1<\cdots<\tau(\ell m)}
 \varepsilon(\tau)\eta^{(1)}_{\{\tau(1),\dots,\tau(m)\}}\cdots\eta^{(k)}_{\{\tau((k-1)m+1),\dots,\tau(km)\}}
 \]
 And so $\left(\Omega_m\left(\mathcal A^{mk}\right)\right)^k=
 \sum_{c\in LQ(m,k)}\varepsilon(c)\eta^{(1)}_{\{1,\dots,mk\}}\cdots\eta^{(k)}_{\{1,\dots,mk\}}$. 
 Proposition \ref{Minors} allows us to conclude.
 $\Box$

 This proposition generalizes a result due to Gherardelli \cite{Ghe} relying the hyperdeterminant of an antisymmetric tensor and the Alon-Tarsi constant \cite{AT}. 
 More precisely, since $LQ(1,k)$ is the set of $k\times k$-latin squares, we recover it for $m=1$,
 $
 \Det\left(\mathcal A^{k}\right)=\frac1{k!}\sum_{c\in LQ(1,k)}\varepsilon(c).$
The following table contains the first values of  $\Pf m\left(\mathcal A^{mk}\right)$.
\[
\begin{array}{|c|cccccccc|}
\hline m\setminus k& 1&2    &3      &4    &5    &6    &7    &8  \\\hline
                  1& 1&1    &0      &4    &0&2304&0&6210846720\\
                  2& 1&3    &90&204120& ?   &?&?&?\\
                  3&0&10&0&?&0&?&0&?\\
                  4&1&35&519750&?&?&?&?&?\\
                  5&0&126&0&?&0&?&0&?\\
                  6&1&462&?&?&?&?&?&?\\\hline  
                  \end{array}
\]
Observe that the first line is the Alon-Tarsi constant for the even values of $k$. The second column is $\binom{2m-1}m$. This can be easily shown by remarking 
that any quasisquare in $LQ(m,2)$ is on the form $$\left[\begin{array}{cccccc}\sigma_1(1)&\cdots&\sigma_1(m)&\sigma_2(1)&\cdots&\sigma_2(m)\\
\sigma_1(m+1)&\cdots&\sigma_1(2m)&\sigma_2(m+1)&\cdots&\sigma_k(2m)\end{array}
\right].$$
We deduce that $\mathrm{card}(LQ(m,2))= \binom{2m}m$ and a straightforward examination shows that any quasisquare has a positive sign.

 When $m$ is even, Theorem \ref{Laplace} allows us to write
\begin{equation}
\Pf m\left(\mathcal A^{mk}\right)=\sum^{\wedge} (-1)^{\Lbag I\Rbag}\Pf m\left(\mathcal A^{mk}[I]\right)\Pf m\left(\mathcal A^{mk}[J]\right),
\end{equation}
where $\mathrm{card}(I^{(s)})=\mathrm{card}(J^{(s)})=\frac m2$ for any $1\leq s\leq k$. Numerical evidences suggest that
 each term having a non-zero contribution in the sum satisfies
$ (-1)^{\Lbag I\Rbag}\Pf m\left(\mathcal A^{mk}[I]\right)\Pf m\left(\mathcal A^{mk}[J]\right)=\Pf m\left(\mathcal A^{mk}[I]\right)^2$. 
If we assume this conjecture, then we show
that $\frac1{k!}\sum_{c\in LQ(m,k)}\varepsilon(c)\geq 0.$ This is still an open problem. For $m=1$, we recover  a weak version of 
the Alon-Tarsi conjecture as stated in \cite{Zappa}.
\section{Concluding remarks}
The construction proposed in this paper allows us to place the Alon-Tarsi conjecture in a broader context. Indeed, a sound knowledge of 
the algebraic
dependences of the different hyperpfaffians for antisymmetric tensors could help us to understand the combinatoric of the Alon-Tarsi sum. The first
(and well known) example is given by $\Pfa=\det^2$ for antisymmetric matrices. This is no longer the case for higher tensors and the complete picture
remains to be discovered.\\
We notice also that there are unsigned version of most of the equalities stated in the paper. These equalities involve hyperhafnians, 
 unsigned analogues of  hyperpfaffians obtained by replacing
 the Grassmann variables by 
commuting nipotent (i.e. $x_i^2=0$) variables. 

 \end{document}